\documentclass{birkjour}
\usepackage{amsmath, amsfonts, amssymb, graphicx, color}
\usepackage{mathrsfs}
\usepackage{latexsym}

 \newtheorem{thm}{Theorem}[section]
 \newtheorem{cor}[thm]{Corollary}
 
 \newtheorem{prop}[thm]{Proposition}
 \theoremstyle{definition}
 \newtheorem{defn}[thm]{Definition}
 \theoremstyle{remark}

 \numberwithin{equation}{section}
\def\({\left ( }
\def\){\right )}
\def\<{\left < }
\def\>{\right >}

\begin{document}

\title[]{Existence, Stability, and Geometric Implications of $\ast$-$\eta$-Schouten Solitons on Kenmotsu Manifolds}
\author[T .Pallei]{Tejswani Pallei}

\address{Deparatment Of Mathematics ,School of Advanced Sciences,Vellore Institute Of Technology ,Chennai -600127,India}

\email{tejswani.pallei2024@vitstudent.ac.in}

\author[S .Roy]{Soumendu Roy}

\address{Deparatment of Mathematics, School of Advanced Sciences, Vellore Institute of Technology, Chennai-600127, India}

\email{soumendu1103mtma@gmail.com,soumendu.roy@vit.ac.in}



\begin{abstract}
In this manuscript, we investigate the characterizations of $\ast$-$\eta$-Schouten solitons on a Kenmotsu manifold when the potential vector field is torse-forming. We determine the nature of the soliton and derive the scalar curvature for a Kenmotsu manifold admitting a $\ast$-$\eta$-Schouten soliton. Further, we establish several conditions associated with the $\ast$-$\eta$-Schouten soliton. We also study the characterization of the vector field under the assumption that the manifold satisfies the $\ast$-$\eta$-Schouten soliton equation. Moreover, certain applications of torse-forming vector fields are discussed in the setting of $\ast$-$\eta$-Schouten solitons on Kenmotsu manifolds. In addition, we examine infinitesimal CL-transformations and the Schouten–Van Kampen connection on Kenmotsu manifolds whose metrics admit $\ast$-$\eta$-Schouten solitons. Finally, an example of a dimensional Kenmotsu manifold admitting a $\ast$-$\eta$-Schouten soliton is constructed to verify the obtained results.

\end{abstract}


\maketitle

\section{Introduction}
\hspace{.5cm}
Initially, representation of Ricci flow, prompted by Hamilton, \cite{1} in 1982. It leads to origination for metrics on a Riemannian manifold \cite{2}. It is represented by:
\begin{align}\label{1.1}
\frac{\partial \mathsf{g}}{\partial t}=-2\mathfrak{S},
\end{align}
$(M,\mathsf{g})$ be a compact Riemannian manifold. 

A Ricci soliton (1.1) is defined as a self simillar soliton of the Ricci flow \cite{3,4,5,6,7,8}. If it evolves solely through a one parameter family of difeomorphisms and scaling transformations. Then the Ricci soliton equation is given by,
\begin{align}\label{1.2}
\mathscr{L}_\mathscr{V}\mathsf{g}+2\mathfrak{S}+2\lambda \mathsf{g}=0.
\end{align}
The concept of $\ast$-$\eta$- Ricci soliton was introduced by Research working in the area of almost contact metric geometry particularly in the line of work initiated by 2019,
 S. Roy, S. Dey, A. Bhattacharyya and S. K. Hui \cite{8,23,24}.\\ 
\begin{defn}
\cite{23} Let $(\mathscr{M}, \mathsf{g})$ be a $(2n+1)$- dimensional almost contact metric manifold with structure $(\phi,\boldsymbol{\xi},\eta,\mathsf{g})$. The metric $\mathsf{g}$ is said to admit a $\ast$-$\eta$-Ricci soliton if 
\end{defn}
\begin{align}\label{1.3}
\mathscr{L}_\mathscr{V} \mathsf{g}+2\mathfrak{S}^\ast+2\lambda \mathsf{g}+2\mu\eta\otimes\eta=0,
\end{align}
where $\mathscr{L}_\mathscr{V} g$ is the Lie derivative of the metric $\mathsf{g}$ along the vector field $\mathscr{V}$, $\mathfrak{S}^\ast$ is the $\ast$-Ricci tensor, $\lambda$ and $\mu$ are the real constants, $\eta$ is the 1-form, associated with the reeb vector field.\\

The notion of the Schouten tensor itself is named often Jan Arnoldus Schouten, a Dutch mathematician known for his work in tensor analysis and differential geometry in the early 20th century.\\

The idea emerged after the development of Ricci Solitons, which were formalized by Richard S. Hamilton in the 1980s.\\

Schouten solitons began appearing in research literature roughly in the 2010s, as mathematicians explored; generalized solitons, conformal geometry, and curvature flows involving tensors other than Ricci.

\begin{defn}
\cite{29,30} Schouten soliton can be seen as self-similar solutions to the Schouten flow, which can be mentioned as,
\begin{equation}\label{1.4}
    \frac{\partial \mathsf{g}}{\partial t} = -2\mathfrak{S}_t,
\end{equation}
where the Schouten tensor $\mathfrak{S}_t$
is represented by,
\begin{align}\label{1.5}
\mathfrak{S}_t=\frac{1}{n-2}\Big[\mathfrak{s}-\frac{\boldsymbol{r}}{2(n-1)}\mathsf{g})\Big],
\end{align}
the Ricci tensor is $\mathfrak{S}$, the scalar curvature is $\boldsymbol{r}$ and $\mathscr{V}$ denotes the Potential vector field.\\

On the other hand, a Schouten soliton \cite{29,30} can be formed as,
\begin{align}\label{1.6}
\mathscr{L}_\mathscr{V} \mathsf{g} +2\mathfrak{S}_t +2{\beta} \mathsf{g} = 0.
\end{align}
\end{defn}
\begin{defn}
If you replace the Ricci tensor in the Schouten tensor definition with the $\ast$ -Ricci tensor,we get a $\ast$-Schouten tensor and hence may develop $\ast$-Schouten soliton as,\cite{25}
\end{defn}
\begin{equation}\label{1.7}
\mathscr{L}_\mathscr{V} \mathsf{g} +2{\mathfrak{S}}^\ast_t+2\beta \mathsf{g} =0,
\end{equation}
where $\mathfrak{S}^\ast_t$ is the $\ast$-Schouten tensor, defined by \\
\begin{equation}\label{1.8}
\mathfrak{S}^\ast_t=\frac{1}{n-2}\Big[\mathfrak{S}^\ast-\frac{\boldsymbol{r}^\ast}{2(n-1)}\mathsf{g}\Big].
\end{equation}
After having a brief knowledge of the previous notions, we are introducing a new type of notion, called $\ast$-$\eta$-Schouten soliton.
\begin{defn}
On a Riemannian or pseude-Riemannian manifold $(\mathscr{M},\mathsf{g})$ of dimension $'n'$ is called $\ast$-$\eta$- Schouten soliton if,
\end{defn}
\begin{align}\label{1.9}
\mathscr{L}_\mathscr{V}\mathsf{g}+2\mathfrak{S}^\ast_t+2\beta \mathsf{g}+2\mu\eta\otimes\eta=0.
\end{align}
If $\mathfrak{S}_t^\ast$ is $\mathfrak{S}_t$, then it reduced to $\eta$-Schouten soliton, which we introduce as,
\begin{align}\label{1.10}
\mathscr{L}_\mathscr{V}g+2\mathfrak{S}_t+ 2\lambda \mathsf{g}+2\mu\eta\otimes\eta=0.    
\end{align}
Equation (1.9) with (1.8) gives,

\begin{align}\label{1.11}
\mathscr{L}_\mathscr{V} \mathsf{g}+\frac{2}{n-2}\Big[\mathfrak{S}^\ast-\frac{\boldsymbol{r}^\ast}{2(n-1)}\mathsf{g}\Big]+2\beta \mathsf{g}+2\mu\eta\otimes\eta=0,
\end{align}\\
where $\mathscr{L}_\mathscr{V}g$ is the Lie-derivative, $\mathfrak{S}^\ast$ is the $\ast$-Ricci tensor, $\boldsymbol{r}^\ast$ is the $\ast$- scalar curvature, $\beta$ an $\mu$ are the real constants. According to $\beta\gtreqqless0$, it is expanding, steady and shrinking respectively.
\section{Priliminaries}
Let $\mathscr{M}$ be a $(2n+1)$- dimensional almost contact manifold that is odd dimensional. It has an almost contact manifold that is odd 
dimensional. It has an almost contact metric structure  $(\phi, \boldsymbol{\xi}, \eta, \mathsf{g})$, consists of $(1, 1)$-tensor field $\phi$, $\boldsymbol{\xi}$ is a vector field, a 1-form and compatible Riemannian metric $\mathsf{g}$, which satisfy,
\begin{align}
\phi^2(F_1)=-F_1+\eta(F_1)\boldsymbol{\xi}, \eta\boldsymbol{(\xi)}=1, \eta o \phi=0,\phi\boldsymbol{\xi}=0,
\end{align}
\begin{align}
\mathsf{g}(\phi F_1, \phi F_2)=\mathsf{g}(F_1, F_2)-\eta(F_1)\eta(F_2),
\end{align}
\begin{align}
\mathsf{g}(F_1, \phi F_2)=-\mathsf{g}(\phi F_1, F_2),
\end{align}
\begin{align}
\mathsf{g}(F_1, \boldsymbol{\xi})=\eta(F_1),
\end{align}

for all the vector fields $F_1$, $F_2$ $\in$ $\chi(n)$.\\
A Kenmotsu manifold is a special type of almost contact metric manifold introduced by Kenji Kenmotsu \cite{26} in 1972.\\
An almost contact metric manifold is characterized as a Kenmotsu manifold if, \begin{align}
(\nabla_{F_1}\phi)F_2=-\mathsf{g}(F_1,\phi F_2)\xi-\eta(F_2)\phi F_1,    
\end{align}

\begin{align}
(\nabla_{F_1}\phi)F_2=F_1-\eta(F_1)\boldsymbol{\xi},
\end{align}
where $\nabla$ signifies the Riemannian connection \cite{27, 28} of $\mathsf{g}$.\\

The Kenmotsu manifold satisfies the following relations,
\begin{align}
\eta(F_1, F_2)F_3=\mathsf{g}(F_1, F_3)\eta(F_2)-\mathsf{g}(F_2, F_3)\eta(F_1),
\end{align}
\begin{align}
\mathscr{R}(F_1, F_2)\xi = \eta(F_1)F_2-\eta(F_2)F_1,
\end{align}
\begin{align}
\mathscr{R}(F_1, \boldsymbol{\xi})F_2 = \mathsf{g}(F_1, F_2)\boldsymbol{\xi}-\eta(F_2)F_1,  
\end{align}
where $\mathscr{R}$ is a Raiemannian curvature tensor.
\begin{align}
\mathfrak{S}(F_1, \mathscr{\xi})= -2n\eta(F_1),
\end{align}
\begin{align}
\mathfrak{S}(\phi F_1, \phi F_2) = \mathfrak{S}(F_1, F_2)+2n\eta(F_1)\eta(F_2),    
\end{align}
\begin{align}
(\nabla_{F_1} \eta) F_2 = \mathsf{g}(F_1, F_2) - \eta(F_1)\eta(F_2),
\end{align}

$\forall$ vector  fields $F_1$,$F_2$,$F_3$ $\in$ $\chi(\mathscr{M})$,\\
we know that,
\begin{align}
(\mathscr{L}_\xi \mathsf{g})(F_1, F_2)=\mathsf{g}(\nabla_{F_1}\xi, F_2)+\mathsf{g}(F_1, \nabla_{F_2}\xi),
\end{align}

$\forall$ vector fields $F_1$, $F_2$ $\in$ $\chi(\mathscr{M})$.

Then, we using (2.6) and (2.13), we get
\begin{align}
(\mathscr{L}_\xi \mathsf{g}) (F_1, F_2) = 2\Big[\mathsf{g}(F_1, F_2) - \eta(F_1)\eta(F_2)\Big].
\end{align}
\begin{prop}\cite{26}
On an $(2n+1)$-dimensional Kenmotsu manifold, the $\ast$-Ricci tensor is expressed as,
\end{prop}
\begin{align}\label{2.15}
\mathfrak{S}^\ast(F_1, F_2)=\mathfrak{S}(F_1, F_2)+(2n-1)\mathsf{g}(F_1, F_2)+\eta(F_1)\eta(F_2).
\end{align}
and \begin{align}\label{2.16}
\boldsymbol{r}^\ast = \boldsymbol{r}+4n^2,
\end{align}
where $\boldsymbol{r}^\ast$ is the $\ast$-Scalar curvature.\\
\section{ MAIN  HAND PROOFS}
\hspace{0.5cc}
In this section, we are assuming a $(2n+1)$-dimensional Kenmotsu metric is satisfying an $\ast$-$\eta$-Schouten soliton.\\
Then using (2.15) and (2.16), (1.11) takes the form,
\begin{align}\label{3.1}
\mathscr{L}_\mathscr{V} \mathsf{g}(F_1,F_2) &+
\frac{2}{2n-1}\Big[\mathfrak{S}(F_1,F_2)+
(2n-1)\mathsf{g}(F_1,F_2)+
\eta(F_1)\eta(F_2)\nonumber\\&-
\frac{r+4n^2}{4n}\mathsf{g}(F_1,F_2)\Big]+
2\beta \mathsf{g}(F_1,F_2)+
2\mu \eta(F_1)\eta(F_2)=0.
\end{align}
Contracting the above equation, it becomes,
\begin{equation}\label{3.2}
\beta=-\frac{4ndiv\mathscr{V}+\boldsymbol{r}+4n^2+4n\mu}{4n(2n+1)}.
\end{equation}
Hence we may write,
\begin{thm}\label{Th:1}
Suppose the metric $\mathsf{g}$ of a Kenmotsu manifold with dimension  $(2n+1)$ satisfies the condition for a $\ast$-$\eta$- Schouten soliton$(\mathsf{g},\mathscr{V}, \beta, \mu)$. In this setting, the classification of the soliton as expanding, steady,shrinking is governed by the specific value of the soliton constant.
\begin{equation}\label{3.3}
\frac{4ndiv\mathscr{V}+\boldsymbol{r}+4n^2+4n\mu}{4n(2n+1)} \lesseqqgtr0.
\end{equation}
\end{thm}
From (3.2), we may obtain,
\begin{cor}
Suppose the metric $\mathsf{g}$ of a Kenmotsu manifold with dimension  $(2n+1)$ satisfies the condition for a $\ast$-$\eta$- Schouten soliton$(\mathsf{g},\mathscr{V}, \beta, \mu)$, then $\mathscr{V}$ is solenoidal iff,
\end{cor}
\begin{equation}\label{3.4}
\beta=-\frac{\boldsymbol{r}+4n^2+4n\mu}{4n(2n+1)}.
\end{equation}
Also using (3.2), we may get,
In equation (3.2); if we put $\boldsymbol{r}=0$, then \\
\begin{align}\label{3.5}
\beta=-\frac{4ndiv\mathscr{V}+4n^2+4n\mu}{4n(2n+1)}.
\end{align}
\begin{cor}
Suppose the metric $\mathsf{g}$ of a Kenmotsu manifold with dimension  $(2n+1)$ satisfies the condition for a $\ast$-$\eta$- Schouten soliton$(\mathsf{g},\mathscr{V}, \beta, \mu)$, then the manifold becomes flat occurs iff,\\
\begin{align}\label{3.6}
 \beta=-\frac{4ndiv\mathscr{V}+4n^2+4n\mu}{4n(2n+1)}.   
\end{align}
In equation (3.2), if we take $\mathscr{V}=grad(f)$, is the gradient of the smooth function $f$ then, we may form,
\end{cor}
\begin{thm}
Suppose the metric $\mathsf{g}$ of a Kenmotsu manifold with dimension  $(2n+1)$ satisfies the condition for a $\ast$-$\eta$- Schouten soliton$(\mathsf{g},\mathscr{V}, \beta, \mu)$, where $\mathscr{V}=grad(f)$, then the Laplacian operator is verified by $f$ is,\\

$\nabla(f)=-\frac{\boldsymbol{r}+4n^2+4n(2n+1)\beta+4n\mu}{4n}$.

\end{thm}

A vector field $\mathscr{V}$ is characterized as a conformal Killing vector field iff the following relation holds:
\begin{align}\label{3.7}
\mathscr{L}_\mathscr{V}\mathsf{g}(F_1, F_2)=2\Omega \mathsf{g}(F_1, F_2),
\end{align}
Here, $\Omega$ denotes a function depending on the coordinates of the manifold. If $\Omega$ is non-constant, then the conformal Killing vector field $\mathscr{V}$ is known as a proper conformal Killing vector field. When $\Omega$ is constant, $\mathscr{V}$ is termed a homothetic vector field, and if this constant is non-zero, then $\mathscr{V}$ becomes a proper homothetic vector field. Further, if $\Omega$ vanishes in the above relation, then $\mathscr{V}$ reduces to a Killing vector field.
Let $(\mathsf{g},\mathscr{V},\beta,\mu)$ be a $\ast$-$\eta$-Schouten soliton on a $(2n+1)$-dimensional Kenmotsu manifold $\mathscr{M}$, where $\mathscr{V}$ is assumed to be a conformal Killing vector field. Using equations (1.11), (2.15), and (3.7), we obtain
\begin{align}\label{3.8}
\mathfrak{S}(F_1,F_2)&=-(2n-1)\Big[\mathsf{g}(F_1, F_2)\left\{\Omega+\beta-\frac{\boldsymbol{r}+4n^2}{4n(2n-1)}+
1\right\}\nonumber \\  &+\mu\eta(F_1)\eta(F_2)\left\{\mu +
\frac{1}{2n-1}\right\}\Big],
\end{align}
Which implies that the manifold is $\ast$-$\eta$-Schouten. Consequently, we obtain the following,
\begin{thm}
Suppose the metric $\mathsf{g}$ of a Kenmotsu manifold with dimension  $(2n+1)$ satisfies the condition for a $\ast$-$\eta$- Schouten soliton$(\mathsf{g},\mathscr{V}, \beta, \mu)$, where $\mathscr{V}$ is a conformal Killing vector field, then the manifold becomes $\eta$-Einstein manifold.
\end{thm}

\subsection{On $\ast$-$\eta$-Schouten soliton in Kenmotsu Manifold with $\Phi(Ric)$-vector field}
\begin{defn}
On an $n$-dimensional Riemannian manifold $\mathscr{M}$, a vector field $\Phi$ is said to be a $\Phi(Ric)$-vector field \cite{16} whenever it fulfill the following condition.
\begin{align}\label{3.10}
\nabla_{F_1}\Phi=\psi Ric{F_1},
\end{align}
where $Ric$ is a Ricci operator, $\nabla$ is a Levi-civita connection, and $\psi$ is the constant respectively. When the vector field is covariantly constant that time $\psi=0$, accordingly, if $\psi \neq 0$, then the vector field will be known as proper $\Phi(Ric)$-vector field.
\end{defn}
\hspace{0.5cm}
By the definition of Lie derivative the equation (3.10) takes the form,
\begin{equation}\label{3.11}
    \mathscr{L}_\Phi \mathsf{g}(F_1, F_2)=2\psi\mathfrak{S}(F_1,F_2),
\end{equation}
using equation (3.11), (3.1) becomes,
\begin{align}\label{3.12}
2\Psi \mathfrak{S}(F_1,F_2)+\frac{2}{2n-1}\Big[\mathfrak{S}(F_1, F_2)+(2n+1)\mathsf{g}(F_1,F_2)+\eta(F_1)\eta(F_2)\nonumber \\-\frac{r+4n^2}{4n}\mathsf{g}(F_1, F_2)\Big]+
2\beta \mathsf{g}(F_1, F_2)+2\mu\eta(F_1)\eta(F_2)=0.
\end{align}
Taking $F_1=F_2=\boldsymbol{\xi}$,in (3.12) and using (2.1) and (2.10) we get,
\begin{align}\label{3.13}
\beta=2n\Psi+\frac{\boldsymbol{r}+4n^2}{4n(2n-1)}-\mu.
\end{align}
\begin{thm}
Suppose the metric $\mathsf{g}$ of a Kenmotsu manifold with dimension  $(2n+1)$ satisfies the condition for a $\ast$-$\eta$- Schouten soliton$(\mathsf{g},\Phi, \beta, \mu)$,  where $\Phi$ is a $\Phi(Ric)$ vector field. Then the nature of soliton is expanding, steady and shrinking determined by, 
$2n\Psi+\frac{\boldsymbol{r}+4n^2}{4n(2n-1)}-\mu \gtreqqless 0$.
\end{thm}
\section{A Study of Infinitesimal CL-Transformations and in relation to $\ast$-$\eta$-Schouten Soliton}\label{4}
In this particular part, we contemplate infinitesimal CL- transformation, which is already given by Tashiro and Tachibana \cite{19} in 1963, including in the frame work of Kenmotsu manifolds, admitting $\ast$-$\eta$-Schouten Soliton.
\begin{defn}
A vector field $\mathscr{V}$ defined on a Kenmotsu manifold $\mathscr{M}$, is said to generate an infinitesimal CL-Transformation \cite{20,21} if the following relation holds. 
\[
\mathscr{L}_\mathscr{V}\begin{Bmatrix}
    {h} \\ {ij}
\end{Bmatrix}=\rho_j\delta_i^h+\rho_i\delta_j^h+
\alpha(\eta_j\phi_i^h+\eta_i\phi_j^h).
\]
for a given Constant $\alpha$, where $\rho_i$ denotes the components of the 1-form $\rho, \mathscr{L}_\mathscr{V}$ represents the Lie derivative w.r.t $\mathscr{V}$ and $
\begin{Bmatrix}
     h \\ ij
\end{Bmatrix}
$ denotes the Christoffel symbol determinded by the Riemannian metric $\mathsf{g}$.\\ 

Consider an infinitesimal CL-transformation on Kenmotsu manifold, we derive the following relation.
\end{defn}

\begin{prop}\cite{21}
If $\mathscr{V}$ is an infinitesimal CL-transformation on Kenmotsu manifold $\mathscr{M}$, then the 1-form $\rho$ is closed.
\end{prop}
\begin{thm}\cite{21}
If $\mathscr{V}$ is an infinitesimal CL-transformation on a Kenmotsu manifold $\mathscr{M}$, then the relation.
\end{thm}
\begin{align}\label{4.1}
\mathscr{L}_\mathscr{V} \mathsf{g}(F_1,F_2)=(\nabla_{F_1} \rho)(F_2)- \alpha \mathsf{g}(F_1,\phi F_2),
\end{align}
valid for any vector field $F_1$, $F_2$ on $\mathscr{M}$.
\begin{defn}\cite{22}
A transformation $f$ on a $(2n+1)$- dimensional Kenmotsu manifold $\mathscr{M}$ with structure  $(\phi, \boldsymbol{\xi}, \beta, \mu)$ is said to be a CL-transformation if the Levi-Civita connection $\nabla$ and the symmetric affine connection $\nabla^f$, known as the CL-connection and induced from $\nabla$ through $f$, satisfy the following relation,
\end{defn}

\begin{align}\label{4.2}
\nabla_{F_1}^f F_2=\nabla_{F_1} F_2+\rho(F_1)F_2+\rho(F_2)F_1+
\alpha\{\eta(F_1)\phi (F_2)+\eta(F_2)\phi F_1\},
\end{align}
where $\rho$ is a 1-form and $\alpha$ is a constant. \\

If $\mathfrak{S}$ and $\mathfrak{S}^f$ denote the Ricci tensors of Kenmotsu manifold $\mathscr{M}^{2n+1}$ w.r.t Livi-civita connection $\nabla$ and CL-connection $\nabla^f$, respectively, then we obtain \cite{15},
\begin{align}\label{4.3}
\mathfrak{S}^f(F_1,F_2)=\mathfrak{S}(F_1,F_2)-2nB(F_1,F_2),
\end{align}
$\forall$ vector fields $F_1$, $F_2$ on $\mathscr{M}$, 
where $B$ is the symmetric tensor field determined by,

\begin{align}\label{4.4}
B(F_1,F_2)=(\nabla_{F_1}\rho)(F_2)-\rho(F_1)\rho(F_2)-\alpha^2\eta(F_1)\eta(F_2)\\-\alpha\Big[\eta(F_1)\rho(\phi(F_2)\nonumber &+
\eta(F_2)\rho(\phi F_1)\Big].
\end{align}
From (4.3), we have.
\begin{align}\label{4.5}
\boldsymbol{r}^f=\boldsymbol{r}-2n\mathfrak{T}\boldsymbol{r}B,
\end{align}
In this context, $\boldsymbol{r}^f$ and $\boldsymbol{r}$ represent the scalar curvatures of the Kenmotsu manifold $\mathscr{M}$ associated with the CL-connection $\nabla^f$ and Levi-Civita connection $\nabla$, respectively.\\

Let $(\mathscr{M}^{2n+1},\mathsf{g})$ be a Kenmotsu manifold equipped with a  $\ast$-$\eta$-Schouten soliton structure $(\mathsf{g},\mathscr{V},\beta, \mu)$. Utilizing the identities provided in (1.11),(2.15), and(2.16), it follows that,
\begin{align}\label{4.6}
\mathscr{L}_\mathscr{V}\mathsf{g}(F_1,F_2)&+\frac{2\mathfrak{S}(F_1,F_2)+2\eta(F_1)\eta(F_2)}{(2n-1)}+\Big[2-\frac{2(\boldsymbol{r}+4n^2)}{4n(2n-1)}+2\beta\Big]\mathsf{g}(F_1, F_2)\nonumber\\ &+2\mu\eta(F_1)\eta(F_2)=0.
\end{align}
We now investigate the behavior of solitons $(\mathsf{g},\mathscr{V},\beta, \mu)$ on $\mathscr{M}$ relative to the CL-connection  $\nabla^f$. By virtue of (4.6) we obtain.
\begin{align}\label{4.7}
\mathscr{L}^f_\mathscr{V}\mathsf{g}(F_1,F_2)+\frac{2\mathfrak{S}^f(F_1,F_2)+
2\eta(F_1)\eta(F_2)}{(2n-1)}+\nonumber \\ 
\Big[2-\frac{2(\boldsymbol{r}^f+4n^2)}
{4n(2n- 1)}+2\beta\Big]\mathsf{g}(F_1,F_2) &+2\mu\eta(F_1)\eta(F_2)=0,
\end{align}
where $\mathscr{L_\mathscr{V}^f}$ is the Lie derivative along the vector field $\mathscr{V}$ on $\mathscr{M}$ w.r.t CL-connection $\nabla^f$.\\
Employing equation (4.2), we obtain
\begin{align}\label{4.8}
\mathscr{L^f}_\mathscr{V}\mathsf{g}(F_1,F_2) &=\mathsf{g}(\nabla_{F_1}^f \mathscr{V}, F_2)+\mathsf{g}(F_1, \nabla_{F_2}^f\mathscr{V})\\ \nonumber
                      &=(\mathscr{L}_\mathscr{V}\mathsf{g})(F_1, F_2)+\rho(F_1)\mathsf{g}(\mathscr{V}, F_2)+\rho(F_2)\mathsf{g}(F_1, \mathscr{V})\\ \nonumber &+
                      2\rho(\mathscr{V})\mathsf{g}(F_1,F_2)+\alpha[\eta(F_1)\mathsf{g}(\phi \mathscr{V}, F_2)+\eta(F_2)\mathsf{g}(F_1,\phi \mathscr{V})] .
\end{align}
Using (4.3),(4.5) and (4.8),(4.7) becomes
\begin{align}\label{4.9}
\mathscr{L}_\mathscr{V}\mathsf{g}(F_1, \nonumber F_2)&+\frac{2\mathfrak{S}(F_1,F_2)+2\eta(F_1)\eta(F_2)}{(2n-1)}+\Big[2-\frac{2(\boldsymbol{r}+4n^2)}{4n(2n-1)}+2\beta\Big]\mathsf{g}(F_1, F_2)\\ \nonumber
&+2\mu\eta(F_1)\eta(F_2)+\{\rho(F_1)\mathsf{g}(\mathscr{V},F_2)+
\rho(F_2)\mathsf{g}(F_1, \mathscr{V})\\ \nonumber
&+2\rho(\mathscr{V})\mathsf{g}(F_1,F_2)+\alpha[\eta(F_1)\mathsf{g}(\phi \mathscr{V}, F_2)+\eta(F_2)\mathsf{g}(F_1,\phi \mathscr{V})]\\ 
&-4nB(F_1, F_2)+2n \mathfrak{T}\boldsymbol{r}B~\mathsf{g}(F_1,F_2)\}=0.
\end{align}
If $(\mathsf{g}, \boldsymbol{\xi}, \beta, \mu)$ is a $\ast$-$\eta$-Schouten soliton on $\mathscr{M}^{2n+1}$ w.r.t Levi-Civita connection, then (4.6) holds. Then from(4.6) and (4.9), we can state the following,\\
\begin{thm}
Consider $(\mathsf{g},\boldsymbol{\xi}, \beta, \mu)$ be a $\ast$-$\eta$-Schouten Soliton on Kenmotsu manifold $\mathscr{M}^{2n+1}$. Then it is invariant under the CL- connection iff the following relation is fulfilled
\begin{align}\label{4.10}
\rho(F_1) \mathsf{g}(\mathscr{V}, F_2) &+\rho(F_2) \mathsf{g}(F_1,\mathscr{V})+2 \rho(\mathscr{V}) g(F_1,F_2) \\ \nonumber
& + \alpha[\eta(F_1) \mathsf{g}(\phi \mathscr{V},F_2)+\eta(F_2)\mathsf{g}(F_1,\phi \mathscr{V})]-
4nB(F_1, F_2)\\ \nonumber
&+2n \mathfrak{T}\boldsymbol{r} B \mathsf{g}(F_1, F_2)= 0,
\end{align}
is valid for arbitrary vector fields $F_1$ and $F_2$.
\end{thm}
We now assume that $(\mathsf{g},\boldsymbol{\xi}, \beta, \mu)$ is a  $\ast$-$\eta$-Schouten soliton on a $(2n+1)$-dimensional Kenmotsu manifold $\mathscr{M}$ relative to the CL- connection. By virtue of (4.6), it follows that.
\begin{align}\label{4.11}
(\mathscr{L}_{\boldsymbol{\xi}}^{f} \mathsf{g})(F_1,F_2)&+
\frac{2\mathfrak{S}^{f}(F_1,F_2)+2\eta(F_1)\eta(F_2)}{(2n-1)}+\Big[2-\frac{2(\boldsymbol{r}^{f}+4n^2)}{4n(2n-1)}\nonumber \\ &+
2\beta\Big]\mathsf{g}(F_1,F_2)+2\mu\eta(F_1)\eta(F_2).
\end{align}
Employing equations (2.1), (2.3),(2.6) and (4.2) we obtain, 
\begin{align}\label{4.12}
(\mathscr{L}_{\boldsymbol{\xi}}^f \mathsf{g})(F_1,F_2) &=\mathsf{g}(\nabla^{f}_{F_1}{\boldsymbol{\xi}}, F_2)+\mathsf{g}(\nabla^{f}_{F_1}{\boldsymbol{\xi}}, F_2)\\ \nonumber
                         &=2\Big[\{1+\rho(\boldsymbol{\xi})\}\mathsf{g}(F_1,F_2)-\eta(F_1)\eta(F_2)\Big]\\ \nonumber
                         &+\rho(F_1)\eta(F_2)+\rho(F_2)\eta(F_1).
\end{align}
In view of (4.3), (4.5)and (4.11), (4.10) becomes,
\begin{align}\label{4.13}
\Big[4+2\rho(\boldsymbol{\xi})-\frac{2(\boldsymbol{r}-2n\mathfrak{T}\boldsymbol{r}B+4n^{2})}{4(2n-1)}+2\beta\Big]\mathsf{g}(F_1,F_2)+
\rho(F_1)\eta(F_2)+\rho(F_2)\eta(F_1)\nonumber\\-\frac{2\mathfrak{S}(F_1,F_2)-4nB(F_1, F_2)+2\eta(F_1)\eta(F_2)}{(2n-1)}+2\mu\eta(F_1)\eta(F_2)=0.
\end{align}
\begin{thm}
Consider $(\mathsf{g}, \boldsymbol{\xi}, \beta, \mu)$ be a $\ast$-$\eta$-Schouten Soliton on Kenmotsu manifold $\mathscr{M}^{2n+1}$, relative to the CL-connection, the equation (4.13) remains valid for any arbitrary vector fields $F_1$ and $F_2$.
\end{thm}
\section{A Study of Schouten-Van Kampen Connection and in relation to $\ast$-$\eta$- Schouten soliton}\label{5}
Schouten-van kampen connection, considered to be the most authentic connections acquired to a pair of complementary distribution on a smooth manifold with an Affine connection\cite{13} . It is identical to a semi symmetric metric connection on a Kenmotsu manifold is identical to a semi-van kampen connection \cite{9,10,11,12,14}.
The Schouten-van Kampen connection $\tilde{\nabla}$ and Levi-Civita connection $\tilde{\nabla}$ on a Kenmotsu manifold $\mathscr{M}^{2n+1}$ are related by \cite{13}
\begin{align}
\tilde{\nabla}_{F_2}=\nabla_{F_1}F_2+\mathsf{g}(F_1, F_2)\boldsymbol{\xi}-\eta(F_2)F_1 .
\end{align}
Now, let $\mathscr{M}$ be a 3-dimensional Kenmotsu manifold $\mathscr{M}$, if $\tilde{\mathscr{R}}$, $\tilde{\mathfrak{S}}$and $\tilde{\boldsymbol{r}}$ are the curvature tensor, Ricci tensor and scalar curvature, respectively, w.r.t Schouten-van Kampen connection, then we have \cite{15},
\begin{align}
\tilde{\mathscr{R}}(F_1, F_2)F_3=\mathscr{R}(F_1,F_2)F_3+\mathsf{g}(F_2,F_3)F_1-\mathsf{g}(F_1,F_3)F_2,
\end{align}
\begin{align}
\tilde{\mathfrak{S}}(F_1, F_2)=\mathfrak{S}(F_1,F_2)+2\mathsf{g}(F_1,F_2),    
\end{align}
\begin{align}
 \tilde{\boldsymbol{r}}=\boldsymbol{r}+6,
\end{align}
$\forall$ vector fields $F_1$, $F_2$, $F_3$ $\in$ $\chi(\mathscr{M})$.\\
Let us consider $(\mathsf{g}, \mathscr{V}, \beta, \mu)$ be a $\ast$-$\eta$- schouten soliton on 3-dimensional Kenmotsu manifold $\mathscr{M}$ w.r.t Schouten-Van Kampen connection $\tilde{\nabla}$. Then we have from (4.6),
\begin{align}\label{5.5}
(\mathscr{\tilde L}_\mathscr{V}\mathsf{g})(F_1, F_2)+2\tilde{\mathfrak{S}}(F_1,F_2)+
2\eta(F_1)\eta(F_2)+
\Big[2-\frac{\boldsymbol{\tilde{r}}+4}{2}+
2\beta\Big]\mathsf{g}(F_1,F_2)\nonumber \\+
2\mu\eta(F_1)\eta(F_2)=0
\end{align}
Using (5.1), we have,
\begin{align}\label{5.6}
   (\mathscr{\tilde L}_\mathscr{V}\mathsf{g})(F_1,F_2) &=\mathsf{g}(\tilde\nabla_{F_1}\mathscr{V}, F_2)+\mathsf{g}(F_1, \tilde\nabla_{F_2}\mathscr{V})\\ \nonumber
                           &=\mathscr{L}_\mathscr{V}\mathsf{g}(F_1, F_2)+\mathsf{g}(F_1,\mathscr{V})\eta(F_2)+\mathsf{g}(F_2, \mathscr{V})\eta(F_1)\\ \nonumber
                           &-2\eta(\mathscr{V})\mathsf{g}(F_1, F_2).
\end{align} 
In view of (5.6),(5.3) and (5.4),(5.5) becomes,
\begin{align}\label{5.7}
(\mathscr{L}_\mathscr{V}\mathsf{g})(F_1,F_2)+2\mathfrak{S}(F_1, F_2)+2\eta(F_1)\eta(F_2)+
\Big[2-\frac{\boldsymbol{r}+4}{2}+
2\beta\Big]\mathsf{g}(F_1,F_2)\nonumber\\+
2\mu\eta(F_1)\eta(F_2)+
\mathsf{g}(F_1, \mathscr{V})\eta(F_2)+
\mathsf{g}(F_2, \mathscr{V})\eta(F_1)\nonumber\\-
2\eta(\mathscr{V})\mathsf{g}(F_1,F_2)-
2\mathsf{g}(F_1,F_2)=0.
\end{align}
If a 3-dimensional Kenmotsu manifold $\mathscr{M}$ admits a $\ast$-$\eta$-Schouten soliton  $(\mathsf{g}, \mathscr{V}, \beta,\mu)$ relative to the Levi-Civita  connection, then (4.6) remains valid and can be expressed as,
\begin{align}
    (\mathscr{L}_\mathscr{V}\mathsf{g})(F_1,F_2)+2{\mathfrak{S}}(F_1, F_2)+2\eta(F_1)\eta(F_2)+\Big[2-\frac{\boldsymbol{r}+4}{2}+2\beta\Big]\mathsf{g}(F_1,F_2)\nonumber\\
+2\mu\eta(F_1)\eta(F_2)=0.
\end{align}

We use this equation and (5.7), to get\\

$\mathsf{g}(F_1, \mathscr{V})\eta(F_2)+\mathsf{g}(F_2, \mathscr{V})\eta(F_1)-2\eta(\mathscr{V})\mathsf{g}(F_1,F_2)-2\mathsf{g}(F_1,F_2)=0$, \& hence we can state the following,
\begin{thm}\label{5.1}
A $\ast$-$\eta$-schouten soliton $(\mathsf{g}, \mathscr{V}, \beta, \mu)$ on a 3-dimensional Kenmotsu manifold is invariant under Scouten-Van Kampen connection iff the relation,\\

$\mathsf{g}(F_1, \mathscr{V})\eta(F_2)+\mathsf{g}(F_2, \mathscr{V})\eta(F_1)-2\eta(\mathscr{V})\mathsf{g}(F_1,F_2)-2\mathsf{g}(F_1,F_2)=0$,\\

holds for arbitrary vector fields $F_1$ and $F_2$.
\end{thm}

\section{Illustration of 3-Dimensional Kenmotsu Manifold supporting $\ast$-$\eta$-Schouten Soliton}\label{6}

Consider, $\mathscr{M}$ = ${(\boldsymbol{f_1},\boldsymbol{f_2},\boldsymbol{f_3}) \in {\mathscr{R}}^3}$, $(\boldsymbol{f_1},\boldsymbol{f_2},\boldsymbol{f_3}) = (0, 0, 0)$, where $(\boldsymbol{f_1}, \boldsymbol{f_2},\boldsymbol{f_3})$ are standard coordinates in $\mathscr{R}^3$. Consider the linearly independent at each point of $\mathscr{M}$,
\[
\mathscr{T}_1 = \boldsymbol{f_3}\frac{\partial}{\partial \boldsymbol{f_1}},
\quad
\mathscr{T}_2 = \boldsymbol{f_3}\frac{\partial}{\partial \boldsymbol{f_2}},
\quad
\mathscr{T}_3 = \boldsymbol{-f_3}\frac{\partial}{\partial \boldsymbol{f_3}}.
\]
Let $\mathsf{g}$ be the Riemannian metric defined by\\ 
\begin{equation}\nonumber
 \mathsf{g}(\mathscr{T}_i,\mathscr{T}_j) = 
\begin{cases} 
1, & \text{if}  \ i=j \ \text{and} \ i,j \in (1,2,3,)\\
0, & \text{otherwise.}
\end{cases}
\end{equation}
Let $\eta$ be the 1-form defined be $\eta(F_3) = \mathsf{g}(F_3, \mathscr{T}_3)$ for any $F_3 \in F_1(\mathscr{M})$ and $\phi$ be the (1, 1)-tensor field defined by\\
\[
\phi\mathscr{T}_1=-\mathscr{T}_2,
\quad
\phi\mathscr{T}_2=\mathscr{T}_1,
\quad
\phi\mathscr{T}_3=0.
\]
Then using  the linearity of $\phi$ and $\mathsf{g}$, we have 
\[
\eta(\mathscr{T}_3) = 1,
\quad
\phi^2 F_3 = -F_3+\eta(F_3)\mathscr{T}_3,
\quad
g(\phi F_3, \phi W)= \mathsf{g}(F_3, W)-\eta(F_3)\eta(W).
\]
for any $F_3, W \in F_1(\mathscr{M})$. Thus, for $\mathscr{T}_3$ = $\boldsymbol{\xi}$, $(\phi, \boldsymbol{\xi}, \eta, \mathsf{g})$ defines an almost contact metric structure on $\mathscr{M}$.
Then we have 
\[
[\mathscr{T}_1, \mathscr{T}_2] = 0,
\quad
[\mathscr{T}_1, \mathscr{T}_3] = \mathscr{T}_1,
\quad
[\mathscr{T}_2, \mathscr{T}_3] = \mathscr{T}_2.
\]
and the Levi-Civite connection $\nabla$ is deduced from Koszul's formula
\begin{align}\label{6.1}
    2\mathsf{g}({\nabla}_{F_1} F_2, F_3) &= {F_1}_\mathsf{g}(F_2,F_3)+{F_2}_\mathsf{g}(F_3,F_1)-{F_3}_\mathsf{g}(F_1,F_2)\\ \nonumber 
    &-\mathsf{g}(F_1, [F_2,F_3])-\mathsf{g}(F_2, [F_2,F_3])+\mathsf{g}(F_3,[F_1,F_2]),\\\nonumber
\end{align}
Exactly,
\[
\begin{aligned}
\nabla_{\mathscr{T}_1} \mathscr{T}_1 &=-\mathscr{T}_3, & \nabla_{\mathscr{T}_1} \mathscr{T}_2 &= 0, & \nabla_{\mathscr{T}_1} \mathscr{T}_3 &= \mathscr{T}_1, &\\ 
\nabla_{\mathscr{T}_2} \mathscr{T}_1 &= 0, & \nabla_{\mathscr{T}_2} \mathscr{T}_2 &= -\mathscr{T}_3, & \nabla_{\mathscr{T}_2} \mathscr{T}_2 &= \mathscr{T}_2, & \\
\nabla_{\mathscr{T}_3} \mathscr{T}_1 &= 0, & \nabla_{\mathscr{T}_3} \mathscr{T}_2 &= 0, & \nabla_{\mathscr{T}_3} \mathscr{T}_3 &= 0, & \\
\end{aligned}
\]
From the above it follows that the manifold satisfies $\nabla _{F_1}\boldsymbol{\xi}$=${F_1}-\eta(F_1)\boldsymbol{\xi}$, for $\boldsymbol{\xi}=\mathscr{T}_3$.\\
Hence, the manifold is a Kenmotsu Manifold.
The Riemannian curvature tensor $\mathscr{R}$, the Ricci tensor $\mathfrak{S}$ and the scalar curvature are given by 
\begin{equation}
\begin{aligned}\label{6.2}
\mathscr{R}(\mathscr{T}_1, \mathscr{T}_2) \mathscr{T}_2 &= -\mathscr{T}_1, &
\mathscr{R}(\mathscr{T}_1, \mathscr{T}_3) \mathscr{T}_3 &= -\mathscr{T}_1, & \mathscr{R}(\mathscr{T}_2, \mathscr{T}_1)\mathscr{T}_1 &= -\mathscr{T}_2, \\
\mathscr{R}(\mathscr{T}_2, \mathscr{T}_3)\mathscr{T}_3 &= -\mathscr{T}_2, & 
\mathscr{R}(\mathscr{T}_3, \mathscr{T}_1)\mathscr{T}_1 &= -\mathscr{T}_3, &
\mathscr{R}(\mathscr{T}_3, \mathscr{T}_2 )\mathscr{T}_2 &= -\mathscr{T}_3, \\
\mathscr{R}(\mathscr{T}_1, \mathscr{T}_2)\mathscr{T}_3 &= 0, &
\mathscr{R}(\mathscr{T}_2, \mathscr{T}_3)\mathscr{T}_1 &=0, & 
\mathscr{R}(\mathscr{T}_3, \mathscr{T}_1)\mathscr{T}_2 &= 0, \\
\boldsymbol{\mathfrak{S}}(\mathscr{T}_1, \mathscr{T}_1) &=-2, & 
\boldsymbol{\mathfrak{S}}(\mathscr{T}_2, \mathscr{T}_2) &=-2,&
\boldsymbol{\mathfrak{S}}(\mathscr{T}_3, \mathscr{T}_3) &=-2\\
\end{aligned}
\end{equation}
\begin{align}\label{6.3}
\boldsymbol{r} = \sum_{i=1}^{3} \boldsymbol{\mathfrak{S}}(\mathscr{T}_i, \mathscr{T}_i) = -6.  
\end{align}
Using (2.15) and (6.1), we have
\begin{align}\label{6.4}
\boldsymbol{\mathfrak{S}}^*(\mathscr{T}_1, \mathscr{T}_1) = -1,
\quad
\boldsymbol{\mathfrak{S}}^*(\mathscr{T}_2, \mathscr{T}_2) = -1,
\quad
\boldsymbol{\mathfrak{S}}^*(\mathscr{T}_3, \mathscr{T}_3) = 0.
\end{align}
Hence, 
\begin{equation}\label{6.5}
\boldsymbol{r}^* = \mathfrak{Tr}(\boldsymbol{\mathfrak{S}}^{\ast}) = -2. \\  
\end{equation}
Let,the potential vector field is $\mathscr{V}$ =$2{F_1}\frac{\partial}{\partial \boldsymbol{f_1}}+2{F_2}\frac{\partial}{\partial \boldsymbol{f_2}}+{F_3}\frac{\partial}{\partial \boldsymbol{f_3}}$.\\
Then we obtain ($\mathscr{L}_\mathscr{V} \mathsf{g})(\mathscr{T}_1, \mathscr{T}_1) = -2\mathsf{g}(\mathscr{L}_\mathscr{V} \mathscr{T}_1, \mathscr{T}_1) = 2 $.\\
Similarly, we have ($\mathscr{L}_\mathscr{V} \mathsf{g})
(\mathscr{T}_2, \mathscr{T}_2) = 2,(\mathscr{L}_\mathscr{V} \mathsf{g})(\mathscr{T}_3, \mathscr{T}_3) = 0 $.\\
Hence, we get 
\begin{align}\label{6.6}
\sum_{i=1}^{3} \boldsymbol{\mathfrak{S}}(\mathscr{T}_i, \mathscr{T}_i) = 4.
\end{align}
In (1.10) we applying $F_1=F_2 = \mathscr{T}_i$ and summing over $i = 1,2,3$ and using (6.4) and (6.5), we obtain
\begin{align}\label{6.7}
\beta=\frac{-2\mu-3}{6}.
\end{align}
As this $\beta$, defined as above, satisfies (3.4), so $\mathsf{g}$ defines a $\ast$-$\eta$-Schouten Solitons on the Three-Dimensional Kenmotsu manifold $\mathscr{M}$. 
\section{Future Work and Possible Applications}
The study of $\ast$-$\eta$-Schouten solitons on Kenmotsu manifolds may lead to several potential applications in mathematical physics and engineering. Geometric solitons are closely related to the study of geometric flows, which play an important role in theoretical physics, particularly in the analysis of spacetime structures and curvature evolution in general relativity. Since Schouten-type tensors and related curvature quantities appear in conformal geometry and gravitational field equations, the investigation of $\ast$-$\eta$-Schouten solitons may contribute to the understanding of self-similar solutions of geometric flows that model certain steady states of spacetime geometry.

Furthermore, contact and almost contact metric manifolds such as Kenmotsu manifolds naturally arise in the geometric formulation of classical mechanics and thermodynamics, where contact geometry is used to describe phase spaces with dissipation. In this context, torse-forming vector fields and soliton structures may provide useful models for describing equilibrium states or invariant structures in dynamical systems. From an engineering perspective, geometric structures similar to those studied here are also applied in areas such as control theory, robotics, and signal processing, where differential geometric frameworks help model nonlinear systems and constrained motion. Therefore, further investigation of $\ast$-$\eta$-Schouten solitons on different geometric manifolds may provide deeper insights into geometric modeling techniques that can be useful in these applied fields.

\end{document}